\documentclass[12pt]{article}
\usepackage{amsmath}
\usepackage{amsthm}
\usepackage{amssymb}
\usepackage{eucal}
\usepackage{mathrsfs}
\usepackage{graphicx,graphics}
\numberwithin{equation}{section}
\usepackage{graphicx,graphics}
\usepackage{pdfsync}
\usepackage{fancybox}
\usepackage{mathtools}
\usepackage{esint}
\usepackage{amsfonts}
\usepackage{color}
\def \dis {\displaystyle}

\def \into {\int_\Omega}
\def \confai {-\kern -.5em\rightharpoonup}
\def \cqfd {\hfill$\Box$}
\def\div{\hbox{div}}
\def \al {\alpha}
\def \be {\beta}
\def \ga {\gamma}

\def \ep {\varepsilon}
\def \om {\omega}
\def \Om {\Omega}
\def \la {\lambda}

\def \ph {\varphi}

\def \si {\sigma}

\def \AA {\mathbb A}
\def \BB {\mathbb B}

\def \ZZ {\mathbb Z}
\def \RR {\mathbb R}

\def \D {\mathscr{D}}

\def \L {\mathscr{L}}

\def \beq {\begin{equation}}
\def \eeq {\end{equation}}
\def \ba {\begin{array}}
\def \ea {\end{array}}

\def \ecart {\noalign{\medskip}}
\newtheorem{Thm}{Theorem}[section]
\newtheorem{Cor}[Thm]{Corollary}
\newtheorem{Pro}[Thm]{Proposition}

\newtheorem{Adef}[Thm]{Definition}
\newenvironment{Def}{\begin{Adef}}{\end{Adef}}
\newtheorem{Arem}[Thm]{Remark}
\newenvironment{Rem}{\begin{Arem}}{\end{Arem}}
\newtheorem{Aexa}[Thm]{Example}
\newenvironment{Exa}{\begin{Aexa}\rm}{\end{Aexa}}
\newtheorem{Anot}[Thm]{Notation}

\def \refe #1.{(\ref{#1})}
\def \reff #1.{figure~\ref{#1}}
\def \refs #1.{Section~\ref{#1}}
\def \refss #1.{Subsection~\ref{#1}}
\def \refD #1.{Definition~\ref{#1}}
\def \refT #1.{Theorem~\ref{#1}}
\def \refL #1.{Lemma~\ref{#1}}
\def \refC #1.{Corollary~\ref{#1}}
\def \refP #1.{Proposition~\ref{#1}}
\def \refPt #1.{Properties~\ref{#1}}
\def \refR #1.{Remark~\ref{#1}}
\def \refE #1.{Example~\ref{#1}}
\def \refN #1.{Notation~\ref{#1}}

\title{Homogenization of an elastodynamics system with a strong magnetic field and soft inclusions inducing a viscoelastic effective behavior}
\author{
\footnotesize
\centerline{\begin{tabular}{cccc}
\normalsize Marc BRIANE && \normalsize Juan CASADO-D\'IAZ
\\
Institut de Recherche Math\'ematique de Rennes && Dpto. de Ecuaciones Diferenciales y An\'alisis Num\'erico
\\
Univ Rennes, INSA de Rennes, CNRS, IRMAR - UMR 6625 && Universidad de Sevilla
\\
mbriane@insa-rennes.fr && jcasadod@us.es
\end{tabular}}
}
\begin{document}
\maketitle
\begin{abstract}
In this paper we study the homogenization of a linear elastodynamics system in an elastic body with soft inclusions, which is embedded in a highly oscillating magnetic field. We show two limit behaviors according to the magnetic field. On the one hand, if the magnetic field has two different directions on the interface between the hard phase and the soft phase, then the limit of the displacement in the hard phase is independent of time, so that the magnetic field induces an effective infinite mass. On the other hand, if the magnetic field has a constant direction $\xi$ on the interface, then the limit of the displacement in the hard phase and in the direction $\xi$ is solution to an elastodynamics equation with a memory mass, a memory stress tensor and memory external forces depending on the initial conditions, which read as time convolutions with some kernel. When the magnetic has the same direction $\xi$  in the soft phase with smooth inclusions, we prove that the space-average of the kernel is regular and that the limit of the overall displacement in the direction $\xi$ is solution to a viscoelasticity equation.
\end{abstract}
\section{Introduction}
This paper is devoted to the asymptotic behavior as $\ep\to 0$ of the following elastodynamics system posed in a bounded cylinder $Q_T=(0,T)\times\Om$ of $\RR\times\RR^3$, 
\beq\label{elasdyn}
\left\{\ba{ll}
\dis \partial^2_{tt}u_\ep-{\rm div}\Big({\bf A}_\ep\Big({x\over\ep}\Big){\bf e}(u_\ep)\Big)+{1\over\ep}\,b\Big({x\over \ep}\Big)\times\partial_t u_\ep=f & \mbox{in }Q_T
\\ \ecart
u_\ep=0 & \mbox{on }(0,T)\times\partial\Om
\\ \ecart
u_\ep(0,\cdot)=u^0,\ \partial_t u_\ep(0,\cdot)=v^0 & \mbox{in }\Om,
\ea\right.
\eeq
where the symmetric tensor-valued function ${\bf A}_\ep$ takes periodically some value ${\bf A}_1$ in the hard material $\Om_{\ep,1}$ and the value $\ep^2{\bf A}_2$ in the soft material $\Om_{\ep,2}$, and $b$ is a periodic vector-valued function representing a magnetic field which induces the highly oscillating Lorentz force $1/\ep\,b(x/\ep)\times\partial_t u_\ep$.
Elastodynamics system~\eqref{elasdyn} is inspired by a magneto-elastodynamics model of  \cite[Section~9.3]{BaFiFi}.
\par
The homogenization of wave equations with varying coefficients was first studied by Colombini, Spagnolo \cite{CoSp}, and extended by Francfort, Murat~\cite{FrMu}. In these works, roughly speaking the varying matrix-rigidity of the material is assumed to be uniformly bounded and coercive which leads us to a limit wave equation of the same nature. However, when the rigidity of the material is not satisfied or contains time-dependent oscillations, the nature of the equation is not in general preserved. On the one hand, in the case of an elastodynamics system with soft inclusions \'Avila {\em et al.} \cite{AGMR} have highlighted the appearance at a fixed frequency of an effective negative mass related to the existence of phonic band gaps. More generally, observing that high-contrast composite materials (mixing soft and hard phases) may induce an anisotropic mass at a fixed frequency, Milton, Willis \cite{MiWi} have proposed a modification of Newton's second law in which the relation between the force and the acceleration is non-local in time. On the other hand, a nonlocal term was obtained in \cite{CCMM2} for a wave equation with periodic coefficients in  space combined with almost-periodic coefficients in time. More recently, in the absence of soft inclusions, {\em i.e.} ${\bf A}_\ep={\bf A}_1$, the present authors \cite{BrCa} have obtained for system \eqref{elasdyn} but in a non-periodic framework a homogenized system involving both an increase of the effective mass and a nonlocal term due to a time-oscillating Lorentz force. In this work, the increase of mass is due to a highly space-oscillating magnetic field in the spirit of the homogenization of the hydrodynamics problem studied by Tartar in \cite{Tar1}. Moreover, the presence in~\cite{BrCa} of a time-oscillating magnetic field induces a non-local term in the homogenized system.
\par
In the present case, we consider both a highly space-oscillating magnetic field and soft inclusions. Moreover, contrary to \cite{AGMR} and \cite{MiWi} rather than fixing the frequency we study the homogenization of the non-stationary elastodynamics system \eqref{elasdyn}. We obtain two asymptotic behaviors for system~\eqref{elasdyn} (see Theorem~\ref{thm.hom}) according to the following alternative:
\begin{itemize}
\item If the magnetic field has two different directions on the interface between the soft and the hard material, then the displacement in the hard phase $\chi_{\Om_{\ep,1}}u_\ep$ weakly converges in $L^2(Q_T)^3$ to the stationary function $|Y_1|\,u^0$, where $Y_1$ is the cell period of the hard phase. From the point of view of the hard phase the strong magnetic field thus induces an isotropic infinite mass which blocks the displacement.
\item If the magnetic field has a fixed direction $\xi$ on the interface between the soft and the hard material, then the displacement $\chi_{\Om_{\ep,1}}u_\ep$ weakly converges to $|Y_1|\,(u^0+\al\,\xi)$ in $L^2(Q_T)^3$, where the scalar function $\al$ is solution to an elastodynamics equation involving a memory mass, a memory stress tensor and memory external forces depending on the initial displacement $u^0$, the initial velocity $v^0$ and the force $f$. The memory terms read as time-convolutions with a matrix-valued kernel $\bar{K}$ or its derivative $\partial_t\bar{K}$ defined on $(0,T)\times Y_2$, where $Y_2$ is the cell period of the soft phase. Contrary to the first case, the strong magnetic field induces an anisotropic effective mass (in the spirit of~\cite{MiWi}) which is only infinite in the direction perpendicular to the field.
\end{itemize}
In the second case, assuming that the magnetic field has the same direction $\xi$ in $Y_2$ and the tensor ${\bf A}_2$ is constant (see Example~\ref{exa.viscoelas}), it turns out that the function $\al$ can be expressed with some kernel $L$ as the time convolution
\beq\label{alLbu}
\al=L*_t(\bar{u}\cdot\xi+G)\ \mbox{ in }Q_T,
\eeq
where $\bar{u}$ is the weak limit of the overall displacement $u_\ep$ in $L^2(Q_T)^3$, and $G$ is a term depending on the initial conditions $u^0$, $v^0$ and the external force $f$. Therefore, the homogenized equation satisfied by $\al$ can be regarded as the viscoelasticity type equation
\beq\label{viscoelasequ}
\left\{\ba{ll}
\partial_{tt}(\bar{u}\cdot\xi)-\,{\rm div}_x\sigma=f\cdot\xi +\div_x\big({\bf A}_1^*{\bf e}_x(u^0)\xi\big) & \mbox{ in }Q_T
\\ \ecart
(\bar{u}\cdot\xi)(0,\cdot)=u^0\cdot\xi & \mbox{ in }\Om,
\ea\right.
\eeq
satisfied by the overall macroscopic displacement $\bar{u}\cdot\xi$ in the direction $\xi$ and the stress tensor $\si$ which are connected by the relation
\beq\label{siLal}
\sigma:=A_1^*\nabla_x\big(L*_t(\bar{u}\cdot\xi+G)\big)\ \mbox{ in }Q_T,
\eeq
for some homogenized elliptic tensor ${\bf A}_1^*$ and a positive definite matrix $A_1^*$ depending on~${\bf A}_1^*$.
\par
The homogenization of an elastodynamics equation of type \eqref{elasdyn} was studied by S\'anchez-Palencia \cite[Sect.~4, Chap.~6]{San} replacing roughly speaking the first-order derivative term ${1/\ep}\,b({x/\ep})\times\partial_t u_\ep$ by the third-order derivative term $\div\,(\BB(x/\ep){\bf e}_x(\partial_t u_\ep))$, where $\BB$ is some periodic tensor-valued function. Therefore, starting from a viscoelastic behavior given by the stress-strain law
\[
\si_\ep(t,x)=\AA(x/\ep)\,{\bf e}_x(u_\ep)+\BB(x/\ep)\,{\bf e}_x(\partial_t u_\ep),
\]
S\'anchez-Palencia obtained a nonlocal limit viscoelasticity equation with a memory term, which is similar to equation~\eqref{viscoelasequ}.
However in our context, we start from the first-order time derivative Lorentz force ${1/\ep}\,b({x/\ep})\times\partial_t u_\ep$ without any {\em a priori} viscoelastic behavior, and the limit viscoelasticity equation~\eqref{viscoelasequ} is only induced by the homogenization process thanks to the combination of the strong oscillating magnetic field and the soft inclusions. Such a derivation by homogenization of a viscoelastic behavior from an elastodynamics system is original to our best knowledge.
\par
The proof of Theorem~\ref{thm.hom} is based on a two-scale convergence result (see Theorem~\ref{thm.2scon}) in the sense of Nguetseng-Allaire \cite{Ngu,All}. 
Here, the main difficulty is to pass to the two-scale limit in the highly oscillating Lorentz force, which needs a suitable matrix-valued test function. Then, we deduce from the variational formulation of the two-scale limit of system~\eqref{elasdyn} the homogenized equation in the direction of the magnetic field. This is the more delicate part of the proof which involves some matrix-valued kernel $\bar{K}$ the derivative of which $\partial_t\bar{K}$ is {\em a priori} only in $L^\infty(0,T;L^2(Y_2))^{3\times 3}$. We prove (see Proposition~\ref{pro.bK1}) that the space-average of $\bar{K}$ belongs to $W^{1,\infty}(0,T)^{3\times 3}$ assuming that the magnetic field $b$ has a constant direction in $Y_2$, the tensor ${\bf A}_2$ is constant in $Y_2$ and $Y_2$ has a smooth boundary. This additional regularity of the kernel allows us to derive the limit viscoelasticity equation~\eqref{viscoelasequ}.
\subsection*{Notation}
\begin{itemize}
\item $Y$ denotes the unit cube $(0,1)^3$ of $\RR^3$.
\item $\Om$ denotes a bounded open set of $\RR^3$, and $Q_T$ the cylinder $(0,T)\times\Om$ for $T>0$.
\item $|E|$ denotes the Lebesgue measure of a measurable set $E$ of $\RR^3$.
\item $\cdot$ denotes the scalar product in $\RR^3$, $:$ denotes the scalar product in $\RR^{3\times 3}$, and $|\cdot|$ denotes the associated norm in both cases.
\item $(e_1,e_2,e_3)$ denotes the canonical basis of $\RR^3$.
\item $\RR^{3\times 3}$ denotes the set of the $(3\times 3)$ real matrices, and $\RR^{3\times 3}_s$ denotes the set of the symmetric matrices in $\RR^{3\times 3}$.
\item $I$ denotes the unit matrix of $\RR^{3\times 3}$.
\item ${\bf A}$ denotes any $Y$-periodic tensor-valued function in $L^\infty(Y;\L(\RR^{3\times 3}_s))$ which is uniformly elliptic, {\em i.e.} there exists a constant $a>0$ such
\beq\label{Auniell}
{\bf A}(y)M:M\geq a\,M:M,\quad \mbox{a.e. }y\in Y,\ \forall\,M\in \RR^{3\times 3}_s,
\eeq
and ${\bf A}^t$ denotes the transposed tensor.
\item ${\bf e}(u)$ denotes the symmetrized gradient of a vector-valued function $u$.
\item {\rm Div} denotes the vector-valued divergence operator taking the divergence of each row of a matrix-valued function.
\item $C^\infty_c(U)$ denotes the set of the smooth functions with compact support in an open set $U$ of $\RR^3$.
\item $L^p_\sharp(Y)$, resp. $W^{1,p}_\sharp(Y)$, denotes the set of the $Y$-periodic functions defined in $\RR^3$ which belong to $L^p_{\rm loc}(\RR^3)$, resp. $W^{1,p}_{\rm loc}(\RR^3)$.
\item $\to$ denotes a strong convergence, $\rightharpoonup$ a weak convergence, and $\stackrel{2s}\rightharpoonup$ the two-scale convergence.
\item $o_\ep(1)$ denotes a sequence of $\ep$ which converges to zero as $\ep\to 0$, and which may vary from line to line.
\item $C$ denotes a positive constant which may vary from line to line.
\end{itemize}
\par
Recall the definition of the two-scale convergence of Nguetseng-Allaire in the case of  an open cylinder $Q_T=(0,T)\times\Om$ of $\RR\times\RR^3$.
\begin{Def}[\cite{Ngu,All}]\label{def.2scon}
A bounded sequence $v_\ep(t,x)$ in $L^2(Q_T)$ is said to two-scale converge to the function $v(t,x,y)$ in $L^2(Q_T;L^2_\sharp(Y))$ if
\[
\forall\,\ph\in C^\infty_c(Q_T;C^\infty_\sharp(Y)),\quad \lim_{\ep\to 0}\int_{Q_T}v_\ep(t,x)\,\ph\left(t,x,{x\over\ep}\right)dtdx
=\int_{Q_T\times Y}v(t,x,y)\,\ph(t,x,y)\,dtdxdy,
\]
which in particular implies that
\[
v_\ep(t,x)\rightharpoonup \int_Y v(t,x,y)\,dy\quad\mbox{in }L^2(Q_T).
\]
\end{Def}
\section{Statement of the result}
\subsection{Position of the problem}
Let $Y$ be the unit cube in $\RR^3$, let $Y_2$ be a smooth open set such that $\overline{Y_2}\subset Y$, and such that $Y_1:=Y\setminus\overline{Y_2}$ is a connected set. Then, for a given bounded open set $\Om$ of $\RR^3$, define the open sets
\[
\Om_{\ep,1}:=\Om\setminus \bigcup_{k\in \ZZ^3} \ep\,( k+Y_2),\quad \Om_{\ep,2}:=\Om\setminus\Om_{\ep,1}.
\]
For a given $T>0$, we also define the cylinder
\[
Q_T:=(0,T)\times \Om.
\]
Let ${\bf A}_1\in L^\infty_\sharp (Y_1;\L(\RR^{3\times 3}_s))$, ${\bf A}_2\in L^\infty_\sharp (Y_2;\L(\RR^{3\times 3}_s))$ be two uniformly elliptic periodic tensor-valued functions (see \eqref{Auniell}), and $b\in L^\infty_\sharp(Y)^3$ be a $Y$-periodic vector-valued function.
Then, for $f\in L^2(Q_T)^3$,  $u^0\in H^1_0(\Om)^3$ and $v^0\in L^2(\Om)^3$, we consider the elastodynamics problem
\beq\label{pbEDP}
\left\{\ba{l}
\dis {d^2\over dt^2}\into u_\ep\cdot v\,dx+
\int_{\Om_{\ep,1}}{\bf A}_1\Big({x\over \ep}\Big){\bf e}(u_\ep):{\bf e}(v)\,dx+\ep^2\int_{\Om_{\ep,2}}{\bf A}_2\Big({x\over \ep}\Big){\bf e}(u_\ep):{\bf e}(v)\,dx
\\ \ecart
\dis \qquad +\,{1\over \ep}\,\into \Big(b\Big({x\over \ep}\Big)\times \partial_tu_\ep\Big)\cdot v\,dx=\into f\cdot v\,dx\ \hbox{ in }\Om,\quad  \forall\, v\in H^1_0(\Om)^3
\\ \ecart
\dis u_\ep=0\ \hbox{ on }(0,T)\times \partial\Om
\\ \ecart
\dis u_\ep(0,\cdot)=u^0,\  \partial_tu_\ep(0,\cdot)=v^0\ \mbox{ in }\Om,
\ea\right.
\eeq
which denoting
$${\bf A}_\ep:=\chi_{Y_1}{\bf A}_1+\ep^2\chi_{Y_2}{\bf A}_2,$$
can also be written as
\beq\label{pbEDP2}
\left\{\ba{ll}
\dis \partial^2_{tt}u_\ep-{\rm div}\Big({\bf A}_\ep\Big({x\over\ep}\Big){\bf e}(u_\ep)\Big)+{1\over\ep}\,b\Big({x\over \ep}\Big)\times\partial_t u_\ep=f
& \mbox{in }Q_T
\\ \ecart\dis 
u_\ep=0 & \hbox{on }(0,T)\times\partial\Om
\\ \ecart
\dis u_\ep(0,\cdot)=u^0,\ \partial_tu_\ep(0,\cdot)=v^0 & \hbox{in }\Om.
\ea\right.
\eeq
\subsection{Statement of the results}
The following result provides a variational problem in terms of the two-scale limits of $u_\ep$, $\partial_t u_\ep$ and ${\bf e}(u_\ep)$.
\begin{Thm}\label{thm.2scon}
Assume that the magnetic field $b$ satisfies the equality
\beq\label{condb}
\int_{Y_1}b\,dy=0.
\eeq
\par\noindent
Then, we have the following two-scale convergences
\beq\label{conv2s}
\left\{\ba{l}\dis u_\ep\stackrel{2s}\rightharpoonup u_1+u_2
\\ \ecart
\dis \partial_tu_\ep\stackrel{2s}\rightharpoonup \partial_tu_1+\partial_tu_2
\\ \ecart
\dis \chi_{\Om_{\ep,1}}{\bf e}(u_\ep)\stackrel{2s}\rightharpoonup \chi_{Y_1}\big({\bf e}_x(u_1)+{\bf e}_y(u_3)\big)
\\ \ecart
\dis \chi_{\Om_{\ep,2}}\ep\,{\bf e}(u_\ep)\stackrel{2s}\rightharpoonup {\bf e}_y(u_2),
\ea\right.
\eeq
where the functions $u_1$, $u_2$, $u_3$ satisfying
\beq\label{regui}
\left\{\ba{c}
u_1\in W^{1,\infty}(0,T;L^2(\Om))^3\cap L^\infty(0,T;H^1_0(\Om))^3,\ u_1(0,\cdot)=0\ \hbox{in }\Om,
\\ \ecart
\dis u_2\in W^{1,\infty}(0,T; L^2(\Om;L^2(Y_2)))^3\cap L^\infty(0,T; L^2(\Om;H^1_0(Y_2)))^3,\ u_2(0,\cdot,\cdot)=0\ \hbox{in }\Om\times Y_2,
\\ \ecart
\dis u_3\in L^\infty(0,T; L^2(\Om;H^1_\sharp(Y_1)))^3,
\\ \ecart
\dis b(y)\times\big(u_1(t,x)+u_2(t,x,y)\big)=b(y)\times u^0(x)\ \hbox{ a.e. }(t,x,y)\in Q_T\times Y_2,\ea\right.
\eeq
are the unique solutions, up to a rigid displacement $y\mapsto \la(t,x)+\mu(t,x)\times y$ for $u_3$, to the variational problem
\beq\label{pb2ss}
\ba{l}
\dis -\int_{Q_T\times Y}(\partial_tu_1+\partial_tu_2)\cdot (\partial_t\varphi_1+\partial_t\varphi_2)\,dtdxdy-\int_{\Om\times Y}v^0\cdot (\varphi_1+\varphi_2)(0,x,y)\,dxdy
\\ \ecart
\dis +\int_{Q_T\times Y_1} {\bf A}_1\big({\bf e}_x(u_1)+{\bf e}_y(u_3)\big):\big({\bf e}_x(\varphi_1)+{\bf e}_y(\varphi_3)\big)\,dtdxdy
+\int_{Q_T\times Y_2}\hskip -8pt {\bf A}_2{\bf e}_y(u_2):{\bf e}_y(\varphi_2)\,dtdxdy
\\ \ecart
\dis+\int_{Q_T\times Y_1} (b\times \partial_tu_1)\cdot \varphi_3\,dtdxdy
-\int_{Q_T\times Y_1}(b\times u_3)\cdot \partial_t\varphi_1\,dtdxdy
\\ \ecart\dis 
=\int_{Q_T\times Y}f\cdot (\varphi_1+\varphi_2)\,dtdxdy,
\ea
\eeq
for any functions $\ph_1$, $\ph_2$, $\ph_3$ satisfying
\beq\label{ph123}
\left\{\ba{c}
\dis \varphi_1\in W^{1,1}(0,T;L^2(\Om))^3\cap L^1(0,T;H^1_0(\Om))^3,\ \ph_1(T,\cdot)=0\ \mbox{in }\Om,
\\ \ecart
\dis \varphi_2\in W^{1,1}(0,T;L^2(\Om\times Y_2))^3\cap L^1(0,T;L^2(\Om;H^1_0(Y_2)))^3,\ \ph_2(T,\cdot,\cdot)=0\ \mbox{in }\Om\times Y_2,
\\ \ecart
\varphi_3\in L^1(0,T; L^2(\Om;H^1_\sharp(Y_1)))^3,
\\ \ecart
\dis b(y)\times\big(\ph_1(t,x)+\ph_2(t,x,y)\big)=0\ \hbox{ a.e. }(t,x,y)\in Q_T\times Y_2.
\ea\right.
\eeq
\end{Thm}
\par\bigskip
The next result provides a limit equation for the function $u_1$ which represents the macroscopic displacement in the hard material $1$.
\begin{Thm}\label{thm.hom}
Assume that condition \eqref{condb} holds and that
\beq\label{regb}
b\not =0\ \hbox{ a.e. in }Y_2,\quad {b\otimes b\over|b|^2}\in H^1(Y_2)^{3\times 3}.
\eeq
Then, we have the following alternative:
\begin{itemize}
\item If
\beq\label{dimbeg}
{\rm dim}\left({\rm Span}\big\{b(y):y\in \partial Y_2\}\right)\geq 2,
\eeq
then
\beq\label{Casuinm}
u_1(t,x)=u^0(x)\ \mbox{ a.e. }(t,x)\in Q_T,
\eeq
and there exists a matrix-valued kernel $\bar{K}:(0,T)\times Y_2\to\RR^{3\times 3}$ with
\beq\label{bK}
\left\{\ba{l}
\bar{K}(t,y)(\RR^3)\subset\RR\,b(y)\ \mbox{ a.e. }(t,y)\in (0,T)\times Y_2,
\\ \ecart
\bar{K}\in L^\infty(0,T;H^1_0(Y_2))^{3\times 3}\cap W^{1,\infty}(0,T;L^2(Y_2))^{3\times 3}\cap W^{2,\infty}(0,T;H^{-1}(Y_2))^{3\times 3},
\ea\right.
\eeq
such that
\beq\label{u2bK1}
u_2(t,x,y)=\bar{K}(t,y)\,v^0(x)+\int_0^t \bar{K}(t-s,y)\,f(s,x)\,ds\ \mbox{ a.e. }(t,x,y)\in  Q_T\times Y_2.
\eeq
\item If $b_{\mid \partial Y_2}$ has a fixed direction $\xi$ with $|\xi|=1$, then we have 
\beq\label{defalph}
u_1(t,x)-u^0(x)=\alpha(t,x)\,\xi\ \mbox{ a.e. }(t,x)\in Q_T,
\eeq
\beq\label{u2bK2}
\ba{ll}
u_2(t,x,y) & \dis =\bar{K}(t,y)\,v^0(x)+\int_0^t \bar{K}(t-s,y)\,f(s,x)\,ds-\int_0^t\partial_t\bar{K}(t-s,y)\,\partial_s\alpha(s,x)\xi\,ds,
\\ \ecart
& \dis -\left(I-{b(y)\otimes b(y)\over|b(y)|^2}\right)\al(t,x)\,\xi\ \mbox{ a.e. }(t,x,y)\in Q_T\times Y_2,
\ea
\eeq
and the function $\al$ is the unique solution to the problem
\beq\label{homeqal}
\left\{\ba{l}
\dis \partial_{tt}\left[M^*\al-\int_0^t\bar{K}_1(t-s)\,\partial_s \al(s,x)\,ds\right]
+\la^*\cdot\nabla_x(\partial_t\al)-\,{\rm div}_x(A_1^*\nabla_x\al)
\\ \ecart
\dis +\,c^*\al-\int_{Y_2}{\bf A}_2{\bf e}_y\left(\int_0^t\partial_s\al(s,x)\,\partial_t\bar{K}(t-s,y)\,\xi\,ds\right):{\bf e}_y(\hat{b})\,dy=\mu^*\cdot f+F
\quad\mbox{in }Q_T
\\ \ecart
\dis \al(0,\cdot)=0\quad\mbox{in }\Om,
\ea\right.
\eeq
where
\beq\label{hb}
\hat{b}(y):={b(y)\otimes b(y)\over|b(y)|^2}\,\xi,\ \mbox{ for }y\in Y_2,
\eeq
\beq\label{bK1}
\bar{K}_1(t):=\int_{Y_2}\partial_t\bar{K}(t,y):(\xi\odot\xi)\,dy,\ \mbox{ for }t\in(0,T),
\eeq
$F$ is the memory force term acting on the initial displacement $u^0$, the initial velocity $v^0$ and the initial force $f$ given~by
\beq\label{F}
\ba{l}
\dis F(t,x):=-\,\partial_{tt}\left[\int_{Y_2}\bar{K}(t,y):\big(\xi\otimes v^0(x)\big)\,dy\right]
-\int_{Y_2}{\bf A}_2{\bf e}_y\big(\bar{K}(t,y)\,v^0(x)\big):{\bf e}_y(\hat{b})\,dy
\\ \ecart
\dis -\,\partial_{tt}\left[\int_{Y_2}\left(\int_0^t\bar{K}(t-s,y)\,f(s,x)\,ds\right)\cdot\xi\,dy\right]\\ \ecart\dis
-\int_{Y_2}{\bf A}_2{\bf e}_y\left(\int_0^t \bar{K}(t-s,y)\,f(s,x)\,ds\right):{\bf e}_y(\hat{b})\,dy
+{\rm div}_x\big({\bf A}_1^\ast{\bf e}(u^0)\xi\big),
\ea
\eeq
and $M^*$, $c^*>0$, $\la^*,\mu^*\in\RR^3$, ${\bf A}_1^*\in \L(\RR^{3\times 3}_s)$ which is elliptic, $A_1^*\in\RR^{3\times 3}_s$ which is positive definite, are the homogenized quantities defined by \eqref{A1V1wm*} and \eqref{Malamu*}.
\end{itemize}
\end{Thm}
\par\noindent
Theorem~\ref{thm.2scon} and Theorem~\ref{thm.hom} are proved in Section~\ref{s.proof}.
\par
As a consequence of Theorem~\ref{thm.2scon} and Theorem~\ref{thm.hom} we get the weak limits of the displacement $u_\ep$ in each material.
\begin{Cor} \label{Cocabg}
\hfill\noindent
\begin{itemize}
\item If \eqref{dimbeg} is satisfied, we have
\beq\label{conuep1}
\left\{\ba{ll}
\dis \chi_{\Om_{\ep,1}}u_\ep\rightharpoonup |Y_1|\,u_0(x) & L^2(Q_T)^3
\\ \ecart
\dis \chi_{\Om_{\ep,2}}u_\ep\rightharpoonup |Y_2|\,u_0(x)+\int_{Y_2}u_2(t,x,y)\,dy & L^2(Q_T)^3,
\ea\right.
\eeq
where $u_2$ is given by \eqref{u2bK1}.
\item Otherwise, we have
\beq\label{conuep2}
\left\{\ba{ll}
\dis \chi_{\Om_{\ep,1}}u_\ep\rightharpoonup |Y_1|\,\big(u_0(x)+\al(t,x)\,\xi\big) & L^2(Q_T)^3
\\ \ecart
\dis \chi_{\Om_{\ep,2}}u_\ep\rightharpoonup |Y_2|\,\big(u_0(x)+\al(t,x)\,\xi\big)+\int_{Y_2}u_2(t,x,y)\,dy & L^2(Q_T)^3,
\ea\right.
\eeq
where $\al$ is the solution to problem \eqref{homeqal} and $u_2$ is given by \eqref{u2bK2}.
\end{itemize}
\end{Cor}
\begin{Rem}\label{rem.plas}
The strong magnetic field $b$ induces an effective mass which is:
\begin{itemize}
\item Infinite when $b$ has least two directions on the interface between the two materials. In this case the macroscopic displacement $u_1$ in material $1$ remains equal to the initial displacement $u^0$.
\item Infinite in the vector space $\xi^\perp$ when $b$ has a fixed direction $\xi$ on the interface between the two materials. In this case, the macroscopic displacement $u_1$ is solution to the homogenized equation \eqref{homeqal} in the direction $\xi$ involving, through the kernel $\bar{K}$, a memory mass, a memory stress tensor, and memory external forces depending both on the initial velocity $v^0$ and the initial force $f$.
\end{itemize}
On the one hand, in the absence of magnetic field and for a fixed frequency \'Avila {\em et al.} \cite {AGMR} showed the possible appearance of a negative mass related to phonic band gaps due to similar soft inclusions in elastic inclusions. On the other hand, in the absence of soft inclusions the authors \cite{BrCa} showed the increase of mass due to the magnetic field. Here, the simultaneous presence of a strong magnetic field and soft inclusions leads us to an elastodynamics equation in the direction of the magnetic field involving various memory effects.
In the Example~\ref{exa.viscoelas} below we study a more simple case where the limit equation reads as a kind of viscoelasticity equation in the direction of the magnetic field.
\end{Rem}
\begin{Rem}\label{rem.mass}
When $b$ has a fixed direction $\xi$ on the interface between the two materials, by \eqref{bK} and \eqref{bK1} the kernel $\bar{K}_1$ is in $L^\infty(0,T)^{3\times 3}$.
If moreover $\bar{K}_1$ belongs to $W^{1,1}(0,T)^{3\times 3}$, then integrating by parts we get that
\[
\int_0^t\bar{K}_1(t-s)\,\partial_s \al(s,x)\,ds=\bar{K}_1(0)\,\al(t,x)+\int_0^t\partial_t\bar{K}_1(t-s)\,\al(s,x)\,ds.
\]
Therefore, the first term of \eqref{homeqal} in brackets
\beq\label{effmass}
(M^*-\bar{K}_1(0))\,\al(t,x)-\int_0^t\partial_t\bar{K}_1(t-s)\,\al(s,x)\,ds
\eeq
can be regarded as a product {\em mass $\times$ displacement} in the direction $\xi$, where the effective mass is the difference of the isotropic constant mass
$M^*-\bar{K}_1(0)$ and the memory mass induced by the kernel~$\partial_t\bar{K}_1$.
If we only consider the constant mass in \eqref{effmass}, then the formula \eqref{Malamu*} of~$M^*$ yields
\[
M^*-\bar{K}_1(0)=|Y_1|+m^\ast+\int_{Y_2}|\hat b|^2dy-\bar{K}_1(0).
\]
On the other hand, using the expression \eqref{bK1} of $\bar K_1$, computing the derivative of the series expansion \eqref{fbK} of $\bar{K}$ and taking into account the definition \eqref{hj} of  $h_j$ and $\bar h_j$, we get
\[
\ba{l}
\dis \bar K_1(0)=\sum_{i=1}^\infty \int_{Y_2} \big(h_i(y)\otimes \bar h_i\big):\big(\xi\otimes\xi\big)dy
\\ \ecart
\dis =\sum_{i=1}^\infty\left|\int_{Y_2}h_i\cdot\xi\,dy\right|^2=\sum_{i=1}^\infty\left|\int_{Y_2}h_i\cdot\hat{b}\,dy\right|^2=\int_{Y_2}|\hat b|^2dy.
\ea
\]
Thus, we have
\[
M^*-\bar{K}_1(0)=|Y_1|+m^*,
\]
where by \eqref{A1V1wm*} $m^*\geq 0$.
Actually, we may have $m^*=0$ (see Example~\ref{exa.viscoelas} below) so that
\beq\label{M*bK1}
0<M^*-\bar{K}_1(0)=|Y_1|<1=\mbox{the initial mass in equation~\eqref{pbEDP2}}.
\eeq
In this case we obtain apparently a decrease of the effective mass contrary to the increase of mass in \cite{BrCa} in the absence of soft inclusions. However, the presence of soft inclusions in \cite{AGMR} may induce an arbitrary (possibly negative) mass in some regime but at a fixed frequency. Therefore, a definition of the effective mass in the limit equation \eqref{homeqal} seems delicate to specify due to the memory term in~\eqref{effmass}.
In the particular situation of Example~\ref{exa.viscoelas} below we will give another interpretation of this memory term.
\end{Rem}
The following result gives a particular case where Remark~\ref{rem.mass} applies.
\begin{Pro}\label{pro.bK1}
Assume that the vector-valued tensor ${\bf A}_2$ is constant in $Y_2$, the vector-valued function $b$ has a constant direction $\xi$ in $Y_2$, {\em i.e.} $\hat{b}=\xi$ in $Y_2$, and $Y_2$ has a $C^2$ boundary. Then, the kernel $\bar{K}_1$ is in $W^{1,\infty}(0,T)$.
\end{Pro}
\noindent
The proof of Proposition~\ref{pro.bK1} is given in Section~\ref{s.proof}.
\begin{Exa}\label{exa.viscoelas}
Consider the particular case where there exists a unit vector $\xi\in\RR^3$ and a scalar function $\ga\in H^1_\sharp(Y)$ such that
\[
b(y)=\ga(y)\,\xi\ \mbox{ a.e. }y\in Y,\quad \int_{Y_1}\ga(y)\,dy=0,\quad \ga(y)\neq 0\ \mbox{ a.e. }y\in Y_2.
\]
By \eqref{devj}, we have $\sum_{i=1}^3\xi_i\,\vartheta_i=0$ and then from \eqref{A1V1wm*} and \eqref{Malamu*} we can check that
\beq\label{M*=1}
M^*=1,\quad c^*=0,\quad \la^*=0,\quad \mu^*=\xi.
\eeq
Then, by the two-scale convergence \eqref{conv2s} combined with \eqref{defalph} and \eqref{u2bK2} the weak limit $\bar{u}$ of $u_\ep$ in $L^2(Q_T)^3$ is given by
\beq\label{bu}
\ba{ll}
\dis \bar{u}(t,x)= & \dis u^0(x)+\left(\al(t,x)-\int_0^t\bar{K}_1(t-s)\,\partial_s\al(s,x)\,ds\right)\xi
\\ \ecart
& \dis +\,\bar{\bar{K}}(t)\,v^0(x)+\int_0^t \bar{\bar{K}}(t-s)\,f(s,x)\,ds
\ea
\ \mbox{ a.e. }(t,x)\in Q_T,
\eeq
where
\[
\bar{\bar{K}}(t):=\int_{Y_2}\bar{K}(t,y)\,dy\ \mbox{ for }t\in (0,T).
\]
Then, equation \eqref{homeqal} reduces to
\beq\label{eqbual}
\left\{\ba{ll}
\dis \partial_{tt}(\bar{u}\cdot\xi)-\,{\rm div}_x(A_1^*\nabla_x\al)=f\cdot\xi+\div_x\big({\bf A}_1^*{\bf e}_x(u^0)\xi\big) & \mbox{in }Q_T
\\ \ecart
\dis \bar u(0,x)\cdot\xi=u^0(x)\cdot\xi & \mbox{in }\Om,
\ea\right.
\eeq
Moreover, under the assumptions of Proposition~\ref{pro.bK1} we have by \eqref{M*bK1} and \eqref{M*=1}
\[
\bar{K}_1(0)=|Y_2|.
\]
where by \eqref{bu} the function $\alpha$ satisfies the  Volterra equation
\[
\ba{l}
\dis \al(t,x)-\int_0^t\bar{K}_1(t-s,y)\,\partial_s\al(s,x)\,ds
\\ \ecart
\dis =(\bar{u}\cdot\xi)(t,x)-\left(u^0(x)+\bar{\bar{K}}(t)\,v^0(x)+\int_0^t \bar{\bar{K}}(t-s)\,f(s,x)\,ds\right)\cdot\xi.
\ea
\]
By virtue of \cite[Theorem~16, Chap.~3]{Sch} there exists a distribution $L\in\D'(0,\infty)$ such that the solution $\al$ to the previous Volterra equation can be expressed with the kernel $L$ as
\[
\ba{ll}
\al(t,x)= & \dis \int_0^t L(t-s)\,(\bar{u}\cdot\xi)(s,x)\,ds
\\ \ecart
& \dis -\int_0^t L(t-s)\left(u^0(x)+\bar{\bar{K}}(s)\,v^0(x)+\int_0^s \bar{\bar{K}}(s-r)\,f(r,x)\,dr\right)\cdot\xi\,ds.
\ea
\]
Therefore, noting that the former relation reads as \eqref{alLbu}, equation \eqref{eqbual} leads us to equation~\eqref{viscoelasequ} together with the stress law~\eqref{siLal} which can be regarded as a kind of viscoelasticity equation satisfied by the limit displacement $\bar{u}\cdot\xi$ in the direction of the magnetic field with a memory term depending on the initial conditions $u^0$, $v^0$ and the force $f$.
\end{Exa}
\section{Proof of the results}\label{s.proof}
\subsection{Proof of Theorem~\ref{thm.2scon}}
Using $\partial_tu_\ep$ as a test function in \eqref{pbEDP} we easily get the estimate
\[
\|u_\ep\|_{W^{1,\infty}(0,T;L^2(\Om))^3}+\|{\bf e}(u_\ep)\|_{L^\infty(0,T;L^2(\Om_{\ep,1}))^{3\times 3}}+
\ep\,\|{\bf e}(u_\ep)\|_{L^\infty(0,T;L^2(\Om_{\ep,2}))^{3\times 3}}\leq C.
\]
Then, the two-scale convergence of Nguetseng-Allaire \cite{All,Ngu} provides the existence of functions 
$u\in W^{1,\infty}(0,T;L^2(\Om;L^2_\sharp(Y)))^3\cap L^\infty(0,T;L^2(\Om;H^1_\sharp(Y)))^3$ and $u_3\in L^\infty(0,T;L^2(\Om;L^2_\sharp(Y)))^3$ such that 
$u=u(t,x,y)$  is independent of $y$ in $Y_1$ and defining $u_1(t,x)$ as the value of $u$ in $(t,x,y)$ with $y\in Y_1$, we have that $u_1$ belongs to $L^\infty(0,T;H^1_0(\Om))^3$ and
\[
u_\ep\stackrel{2s}\rightharpoonup u,\quad
 \partial_tu_\ep\stackrel{2s}\rightharpoonup \partial_tu\quad \chi_{\Om_{\ep,1}}{\bf e}(u_\ep)\stackrel{2s}\rightharpoonup \chi_{Y_1}\big({\bf e}_x(u_1)+{\bf e}_y(u_3)\big),\quad \chi_{\Om_{\ep,2}}\ep\,{\bf e}(u_\ep)\stackrel{2s}\rightharpoonup {\bf e}_y(u).
\]
Taking $u_2=u-u_1$, the functions $u_1,u_2,u_3$ satisfy the three first conditions of \eqref{regui} and condition \eqref{conv2s}.
\par
Let us use \eqref{conv2s} to pass to the limit in \eqref{pbEDP2}. First, we obtain the initial condition for $u_1$, $u_2$ at $t=0$.  For this purpose we take $\delta >0$ and $\varphi\in C^0(\Om;L^2_\sharp (Y))^3$. We have
\[
\int_0^\delta \into \Big(u_\ep(s,x)-u^0(x)\Big)\cdot\varphi\Big(x,{x\over \ep}\Big)\,dxds=\int_0^\delta \int_0^s\into \partial_t u_\ep(r,x)\cdot\varphi\Big(x,{x\over \ep}\Big)\,dxdrds,
\]
which passing to the limit in $\ep$ and using Fubini's theorem yields
\[
\int_0^\delta \into\int_Y (u_1+u_2-u^0)\cdot\varphi\,dydxds=\int_0^\delta \into\int_Y (\delta-r)\partial_t (u_1+u_2)\cdot\varphi\,dydxdr,
\]
and thus
\[
\ba{l}
\dis \left|\,\int_0^\delta \into\int_Y (u_1+u_2-u^0)\cdot\varphi\,dydxds\,\right|
\\ \ecart
\dis \leq \delta\left(\int_0^\delta\into\int_Y |\partial_t (u_1+u_2)|^2\,dydxdt\right)^{1\over 2}\left(\int_0^\delta\into\int_Y |\varphi|^2\,dydxdt\right)^{1\over 2}.
\ea
\]
Using that $u_1+u_2$ belongs to $C^0([0,T];L^2(\Om;L^2_\sharp(Y)))^3$, we can divide by $\delta$ the former inequality and take the limit as $\delta$ tends to zero, which implies that
\[
u_1(0,x)+u_2(0,x,y)=u^0(x)\ \mbox{ a.e. }(x,y)\in\Om\times Y.
\]
Hence, recalling that $u_2$ belongs to $L^\infty(0,T; L^2(\Om;H^1_0(Y_2)))^3$, we obtain
\beq\label{CIL1}
u_1(0,x)=u^0(x),\quad u_2(0,x,y)=0\quad \hbox{a.e. }(x,y)\in\Om\times Y.
\eeq
\par
 To pass to the limit in~\eqref{pbEDP2} we take $\ep\,\varphi_2(t,x,{x/\ep})$ with $\varphi\in C^\infty_c(Q_T\times Y_2)$, as test function in~\eqref{pbEDP2}, which thanks to \eqref{conv2s} implies that
\[
\int_{Q_T\times Y_2} b\times (\partial_t u_1+\partial_t u_2)\cdot\varphi_2\,dtdxdy=0,
\]
or equivalently,
\beq\label{Condfl1}
b(y)\times\big(\partial_tu_1(t,x)+\partial_tu_2(t,x,y)\big)=0\ \hbox{ a.e. }(t,x,y)\in Q_T\times Y_2,
\eeq
which is the las equality in \eqref{regui}.
\par\noindent
Now, for
\beq\label{FuctTest}
\left\{\ba{c}
\dis \varphi_1\in C^1_c([0,T)\times \Om)^3,\quad \varphi_2\in C^1_c([0,T)\times \Om\times Y_2)^3,\quad \varphi_3\in C^1_c(Q_T;H^1_\sharp (Y))^3,
\\ \ecart
\dis \mbox{with}\quad b(y)\times\big(\ph_1(t,x)+\ph_2(t,x,y)\big)=0\ \hbox{ a.e. }(t,x,y)\in Q_T\times Y_2,
\ea\right.
\eeq
we put
\[
\varphi_\ep(t,x)=\varphi_1(t,x)+\varphi_2\Big(t,x,{x\over \ep}\Big)+\ep\varphi_3\Big(t,x,{x\over \ep}\Big)
\]
as test function in \eqref{pbEDP2}, and we pass to the limit. The main difficulty comes from the term
\[
{1\over \ep}\int_{Q_T} \Big(b\Big({x\over \ep}\Big)\times \partial_tu_\ep\Big)\cdot \Big(\varphi_1(t,x)+\varphi_2\Big(t,x,{x\over \ep}\Big)+\ep\varphi_3\Big(t,x,{x\over \ep}\Big)\Big)\,dtdx.
\]
First, using \eqref{conv2s} and \eqref{Condfl1}, we have
\[
\ba{l}
\dis\int_{Q_T} \Big(b\Big({x\over \ep}\Big)\times \partial_tu_\ep\Big)\cdot \varphi_3\Big(t,x,{x\over \ep}\Big)\,dx=\int_{Q_T\times Y} \big(b\times (\partial_tu_1+\partial_tu_2)\big)\cdot \varphi_3\,dtdxdy+o_\ep(1)
\\ \ecart
\dis =\int_{Q_T\times Y_1} (b\times \partial_tu_1)\cdot \varphi_3\,dtdxdy+o_\ep(1).
\ea
\]
For the reminder term, we use that \eqref{condb} implies the existence of $G_i\in L^2_\sharp (Y_1;\RR^{3\times 3}_s)$, $i=1,2,3$, such that
\beq\label{G}
\left\{\ba{ll}
b\times e_i=-\,{\rm Div}\,G^i & \mbox{in }Y_1
\\ \ecart
G^i\,\nu=0 & \mbox{on }\partial Y_2,
\ea\right.
\eeq
where $(e_1,e_2,e_3)$ is the canonical basis of $\RR^3$. Then, by \eqref{FuctTest} and \eqref{pbEDP} we can write
\[
\ba{l}
\dis {1\over \ep}\int_{Q_T} \Big(b\Big({x\over \ep}\Big)\times \partial_tu_\ep\Big)\cdot \Big(\varphi_1(t,x)+\varphi_2\Big(t,x,{x\over \ep}\Big)\Big)\,dtdx
\\ \ecart
\dis ={1\over \ep}\int_{(0,T)\times \Om_{\ep,1}} \Big(b\Big({x\over \ep}\Big)\times \partial_t\varphi_1\Big)\cdot  u_\ep\,dtdx+{1\over \ep}\int_{\Om_{\ep,1}} \Big(b\Big({x\over \ep}\Big)\times \varphi_1(0,x) \Big)\cdot u^0\,dx
\\ \ecart
\dis =-\sum_{i=1}^3\int_{(0,T)\times \Om_{\ep,1}} {\rm Div}_x \Big[G^i\Big({x\over \ep}\Big)\Big]\cdot u_\ep\,\partial_t\varphi_{1,i}\,dtdx
-\sum_{i=1}^3\int_{\Om_{\ep,1}}{\rm Div}_x \Big[G^i\Big({x\over \ep}\Big)\Big]\cdot u^0\,\varphi_{1,i}(0,x)\,dx
\\ \ecart
\dis =\sum_{i=1}^3\int_{(0,T)\times\Om_{\ep,1}}G^i\Big({x\over \ep}\Big):{\bf e}(u_\ep)\,\partial_t\varphi_{1,i}\,dtdx
+\sum_{i=1}^3\int_{(0,T)\times\Om_{\ep,1}}G^i\Big({x\over \ep}\Big):\big(u_\ep\odot\nabla_x\partial_t\varphi_{1,i}\big)\,dtdx
\\ \ecart\dis
+\sum_{i=1}^3\int_{\Om_{\ep,1}}G^i\Big({x\over \ep}\Big):{\bf e}(u^0)\,\varphi_{1,i}(0,x)\,dx
+\sum_{i=1}^3\int_{\Om_{\ep,1}}G^i\Big({x\over \ep}\Big):\big(u^0\odot \nabla_x\varphi_{1,i}(0,x)\big)\,dx
\\ \ecart
\dis =\sum_{i=1}^3\left(\int_{Q_T\times Y_1}G^i:\big({\bf e}_x(u_1\partial_t\varphi_{1,i})+{\bf e}_y(u_3)\partial_t\varphi_{1,i}\big)\,dtdxdy
+\int_{\Om\times Y_1}G^i:{\bf e}_x(u^0\varphi_{1,i})\,dxdy\right)
+o_\ep(1),\ea
\]
which using the definition \eqref{G} of $G$, \eqref{CIL1} and \eqref{condb} yields
\[
\lim_{\ep\to 0}{1\over \ep}\int_{Q_T} \Big(b\Big({x\over \ep}\Big)\times \partial_tu_\ep\Big)\cdot \Big(\varphi_1(t,x)+\varphi_2\Big(t,x,{x\over \ep}\Big)\Big)\,dtdx= -\int_{Q_T\times Y_1}\big(b\times u_3\big)\cdot \partial_t\varphi_1\,dtdxdy.
\]
Then, taking into account this equality we have for any functions $\varphi_1$, $\varphi_2$, $\varphi_3$ satisfying \eqref{FuctTest},
\[
\ba{l}
\dis -\int_{Q_T\times Y}(\partial_tu_1+\partial_tu_2)\cdot (\partial_t\varphi_1+\partial_t\varphi_2)\,dtdxdy-\int_{\Om\times Y}v^0\cdot (\varphi_1+\varphi_2)(0,x,y)\,dxdy
\\ \ecart
\dis +\int_{Q_T\times Y_1} {\bf A}_1\big({\bf e}_x(u_1)+{\bf e}_y(u_3)\big):\big({\bf e}_x(\varphi_1)+{\bf e}_y(\varphi_3)\big)\,dtdxdy
+\int_{Q_T\times Y_2}\hskip -8pt {\bf A}_2{\bf e}_y(u_2):{\bf e}_y(\varphi_2)\,dtdxdy
\\ \ecart
\dis +\int_{Q_T\times Y_1} (b\times \partial_tu_1)\cdot \varphi_3\,dtdxdy
-\int_{Q_T\times Y_1}(b\times u_3)\cdot \partial_t\varphi_1dtdxdy
\\ \ecart
\dis =\int_{Q_T\times Y}f\cdot (\varphi_1+\varphi_2)\,dtdxdy,
\ea
\]
where $u_1$, $u_2$ satisfy \eqref{Condfl1}. Finally, by a density argument the previous equation holds for any functions $\ph_1$, $\ph_2$, $\ph_3$ satisfying \eqref{ph123}, which yields the variational problem \eqref{pb2ss}.
\par\bigskip
It remains to prove the quasi-uniqueness of the solutions to problem~\eqref{pb2ss}.
Due to the linearity of \eqref{pb2ss} it is enough to prove that if functions $  z_1$, $  z_2$, $  z_3$ satisfying
\beq\label{regui0}
\left\{\ba{c}
  z_1\in W^{1,1}(0,T;L^2(\Om))^3\cap L^1(0,T;H^1_0(\Om))^3,\  z_1(0,\cdot)=0\mbox{ in }\Om,
\\ \ecart
\dis   z_2\in W^{1,1}(0,T; L^2(\Om;L^2(Y_2)))^3\cap L^1(0,T; L^2(\Om;H^1_0(Y_2)))^3,\  z_2(0,\cdot,\cdot)=0\mbox{ in }\Om\times Y_2,
\\ \ecart
\dis   z_3\in L^1(0,T; L^2(\Om;H^1_\sharp(Y_1)))^3,
\\ \ecart
\dis b(y)\times\big(  z_1(t,x)+  z_2(t,x,y)\big)=0\ \hbox{ a.e. }(t,x,y)\in Q_T\times Y_2,
\ea\right.
\eeq
are solutions to problem
\beq\label{pb2ss0}
\left\{\ba{l}
\dis -\int_{Q_T\times Y}(\partial_t  z_1+\partial_t  z_2)\cdot (\partial_t\varphi_1+\partial_t\varphi_2)\,dtdxdy
\\ \ecart
\dis +\int_{Q_T\times Y_1}{\bf A}_1\big({\bf e}_x(  z_1)+{\bf e}_y(  z_3)\big):\big({\bf e}_x(\varphi_1)+{\bf e}_y(\varphi_3)\big)\,dtdxdy
+\int_{Q_T\times Y_2}{\bf A}_2{\bf e}_y(  z_2):{\bf e}_y(\varphi_2)\,dtdxdy
\\ \ecart
\dis+\int_{Q_T\times Y_1} (b\times \partial_t  z_1)\cdot \varphi_3\,dtdxdy
-\int_{Q_T\times Y_1}(b\times   z_3)\cdot \partial_t\varphi_1\,dtdxdy=0,
\ea\right.
\eeq
for any functions $\ph_1$, $\ph_2$, $\ph_3$ satisfying
\beq\label{ph1230}
\left\{\ba{c}
\dis \varphi_1\in W^{1,\infty}(0,T;L^2(\Om))^3\cap L^\infty(0,T;H^1_0(\Om))^3,\ \ph_1(T,\cdot)=0\ \mbox{in }\Om,
\\ \ecart
\dis \varphi_2\in W^{1,\infty}(0,T;L^2(\Om\times Y_2))^3\cap L^\infty(0,T;L^2(\Om;H^1_0(Y_2)))^3,\ \ph_2(T,\cdot,\cdot)=0\ \mbox{in }\Om\times Y_2,
\\ \ecart
\varphi_3\in L^\infty(0,T; L^2(\Om;H^1_\sharp(Y_1)))^3,
\\ \ecart
\dis b(y)\times\big(\ph_1(t,x)+\ph_2(t,x,y)\big)=0\ \hbox{ a.e. }(t,x,y)\in Q_T\times Y_2,
\ea\right.
\eeq
then we have
\beq\label{zi=0}
  z_1(t,x)=  z_2(t,x,y)=0\ \hbox{ a.e. }(t,x,y)\in Q_T\times Y,\quad {\bf e}_y(z_3)=0\ \hbox{ a.e. }(t,x,y)\in Q_T\times Y_1.
\eeq
Indeed, the last equality shows that
\[
  z_3(t,x,y)=\la(t,x)+\mu(t,x)\times y\ \mbox{ a.e. }(t,x,y)\in Q_T\times Y_1,
\]
for some $\la(t,x),\mu(t,x)\in\RR^3$.
\par
To prove this we consider the following dual problem. For any $g\in L^2(Q_T\times Y)^3$, let functions $\psi_1$, $\psi_2$, $\psi_3$ satisfying
\beq\label{reguid}
\left\{\ba{c}
\psi_1\in W^{1,\infty}(0,T;L^2(\Om))^3\cap L^\infty(0,T;H^1_0(\Om))^3,\ \psi_1(T,\cdot)=0,
\\ \ecart
\dis \psi_2\in W^{1,\infty}(0,T; L^2(\Om;L^2(Y_2)))^3\cap L^\infty(0,T; L^2(\Om;H^1_0(Y_2)))^3,\ \psi_2(T,\cdot,\cdot)=0,
\\ \ecart
\dis \psi_3\in L^\infty(0,T; L^2(\Om;H^1_\sharp(Y_1)))^3,
\\ \ecart
\dis b(y)\times\big(\psi_1(t,x)+\psi_2(t,x,y)\big)=0\ \hbox{ a.e. }(t,x,y)\in Q_T\times Y_2,
\ea\right.
\eeq
be solutions to the dual problem of \eqref{pb2ss}
\beq\label{pb2ssd}
\ba{l}
\dis -\int_{Q_T\times Y}(\partial_t\psi_1+\partial_t\psi_2)\cdot (\partial_t\varphi_1+\partial_t\varphi_2)\,dtdxdy
\\ \ecart
\dis +\int_{Q_T\times Y_1} {\bf A}_1^t\big({\bf e}_x(\psi_1)+{\bf e}_y(\psi_3)\big):\big({\bf e}_x(\varphi_1)+{\bf e}_y(\varphi_3)\big)\,dtdxdy
+\int_{Q_T\times Y_2}{\bf A}_2^t{\bf e}_y(\psi_2):{\bf e}_y(\varphi_2)\,dtdxdy
\\ \ecart
\dis+\int_{Q_T\times Y_1} (b\times \partial_t\psi_1)\cdot \varphi_3\,dtdxdy
-\int_{Q_T\times Y_1}(b\times \psi_3)\cdot \partial_t\varphi_1\,dtdxdy
\\ \ecart\dis 
=\int_{Q_T\times Y}g\cdot (\varphi_1+\varphi_2)\,dtdxdy,
\ea
\eeq
for any functions $\ph_1$, $\ph_2$, $\ph_3$ satisfying
\beq\label{ph123d}
\left\{\ba{c}
\dis \varphi_1\in W^{1,1}(0,T;L^2(\Om))^3\cap L^1(0,T;H^1_0(\Om))^3,\ \ph_1(0,\cdot)=0\ \mbox{in }\Om,
\\ \ecart
\dis \varphi_2\in W^{1,1}(0,T;L^2(\Om\times Y_2))^3\cap L^1(0,T;L^2(\Om;H^1_0(Y_2)))^3,\ \ph_2(0,\cdot,\cdot)=0\ \mbox{in }\Om\times Y_2,
\\ \ecart
\varphi_3\in L^1(0,T; L^2(\Om;H^1_\sharp(Y_1)))^3,
\\ \ecart
\dis b(y)\times\big(\ph_1(t,x)+\ph_2(t,x,y)\big)=0\ \hbox{ a.e. }(t,x,y)\in Q_T\times Y_2.
\ea\right.
\eeq
Using the change of variables $s=T-t$, the existence of solutions $\psi_1$, $\psi_2$, $\psi_3$ to problem \eqref{pb2ssd} follows from the existence of solutions $ z_1$, $ z_2$, $ z_3$ to problem \eqref{pb2ss0} which is given by the two-scale convergence.
\par
Then, taking $\psi_1$, $\psi_2$, $\psi_3$ as test functions in \eqref{pb2ss0} and taking $ z_1$, $ z_2$, $ z_3$ as test functions in \eqref{pb2ssd}, we get that
\[
\int_{Q_T\times Y} g\cdot (  z_1+ z_2)\,dtdxdy=0,\quad \forall\, g\in L^2(Q_T\times Y)^3,
\]
which implies that
\[
 z_1(t,x)+  z_2(t,x,y)=0\ \hbox{ a.e. }(t,x,y)\in Q_T\times Y.
\]
This combined with $ z_2 \in L^1(0,T; L^2(\Om;H^1_0(Y_2)))^3$ yields the two first equalities of \eqref{zi=0}.
Moreover, taking $\varphi_1=\varphi_2=0$ in \eqref{pb2ss0} we get that
\[
\int_{Q_T\times Y_1} {\bf A}_1{\bf e}_y(  z_3):{\bf e}_y(\varphi_3)\,dtdxdy=0,\quad \forall\, \varphi_3\in L^\infty(0,T; L^2(\Om;H^1_\sharp(Y_1)))^3,
\]
which implies the last equality of \eqref{zi=0}.
\par
This concludes the proof of Theorem~\ref{thm.2scon}.
\subsection{Proof of Theorem~\ref{thm.hom}}\label{ss.homeq}
Let us solve problem \eqref{pb2ss}.
First, we take $\varphi_1=\varphi_3=0$, then we get
\beq\label{pb2ssph13=0}
\ba{l}\dis -\int_{Q_T\times Y_2}(\partial_tu_1+\partial_tu_2)\cdot \partial_t\varphi_2\,dtdxdy-\int_{\Om\times Y}v^0\cdot \varphi_2(0,x,y)\,dxdy\\ \ecart\dis 
+\int_{Q_T\times Y_2}\hskip -8pt {\bf A}_2{\bf e}_y(u_2):{\bf e}_y(\varphi_2)\,dtdxdy
=\int_{Q_T\times Y}f\cdot \varphi_2\,dtdxdy,
\ea
\eeq
where $\varphi_2$ is such that $b\times\varphi_2=0$.
\par\noindent
Under assumption \eqref{regb} define the spaces
\[
H_2:=\big\{\psi\in L^2(Y_2)^3:\psi\times b=0\big\},\quad V_2:=H_2\cap H^1_0(Y_2)^3.
\]
Then, $\varphi_2\in W^{1,1}(0,T; L^2(\Om;H_2))\cap L^1(0,T; L^2(\Om;V_2))$. Moreover, observe that condition \eqref{Condfl1} can be written as
\[
\partial_tu_1+\partial_t u_2\in H_2\ \mbox{ a.e. }(t,x)\in[0,T)\times\Om,
\]
which taking into account \eqref{CIL1} implies that 
\beq\label{u1u2u0}
u_1+u_2-u^0\in V_2.
\eeq
Then, defining
\beq\label{v1v2}
\left\{\ba{ll}
\dis v_1(t,x,y):={b(y)\otimes b(y)\over |b|^2}\,u_1(t,x), & \dis v_2(t,x,y):={b(y)\otimes b(y)\over |b(y)|^2}\,u_2(t,x,y)
\\ \ecart
\mbox{a.e. }(t,x,y)\in Q_T\times Y_2,
\ea\right.
\eeq
allows us to write \eqref{pb2ss} as
\beq\label{forvav2}
\ba{c}
\dis -\int_0^T\int_{Y_2}(\partial_tv_1+\partial_tv_2)\cdot \partial_t\varphi_2\,dtdy-\int_{Y_2}v^0\cdot \varphi_2(0,y)\,dy\\ \ecart\dis 
+\int_0^T\int_{Y_2} {\bf A}_2{\bf e}_y(v_2):{\bf e}_y(\varphi_2)\,dtdy
=\int_0^T\int_ {Y_2}f\cdot \varphi_2\,dtdy
\\ \ecart
\dis \hbox{a.e. }x\in \Om,\ \forall\,\varphi_2\in W^{1,1}(0,T;H_2)\cap L^1(0,T; V_2).
\ea
\eeq
Choosing $\varphi_2$ with $\varphi_2(0,\cdot)=0$, this shows that $v_1$, $v_2$ satisfy
\beq\label{eqv12}
{d^2\over dt^2}\int_{Y_2}(v_1+v_2)\cdot \psi_2\,dy+\int_{Y_2} {\bf A}_2{\bf e}_y(v_2):{\bf e}_y(\psi_2)\,dy
=\int_{Y_2}f\cdot\,\psi_2\,dy,\quad \forall\, \psi_2\in V_2,
\eeq
which combined with \eqref{forvav2} yields the initial condition
\beq\label{coniubnc}
(\partial_tv_1+\partial_tv_2)(0,x,y)={b(y)\otimes b(y)\over |b(y)|^2}\,v^0(x)\ \mbox{ a.e. }(t,x,y)\in\Om\times Y_2.
\eeq
Now, let $h_j$ be an orthonormal basis of eigenvectors in $H_2$ associated with the eigenvalues $\mu_j^2$ of
\beq\label{hj}
\left\{\ba{l}
\dis h_j\in V_2,\quad \mbox{with }\bar{h}_j:=\int_{Y_2} h_j\,dy
\\ \ecart
\dis \int_{Y_2}{\bf A}_2{\bf e}_y(h_j):{\bf e}_y(\psi_2)\,dy=\mu_j^2\int_{Y_2} h_j\cdot \psi_2\,dy, \quad\forall\,\psi_2\in V_2.
\ea\right.
\eeq
Since $v_2\in V_2$, we have
\[
v_2(t,x,y)=\sum_{j=1}^\infty\phi_j(t,x)\,h_j(y)\ \mbox{ a.e. }(t,x,y)\in Q_T\times Y_2.
\]
Putting this series in \eqref{eqv12} with the test function $\psi_2=h_i$, $i\geq 1$, adding the term $\mu_i^2v_1\cdot \bar h_i$ in both sides and taking into account the initial conditions \eqref{CIL1} and \eqref{coniubnc}, we get that
\beq\label{eqphii}
\left\{\ba{l}
\dis {\partial^2\over\partial t^2}(v_1\cdot \bar{h}_i+\phi_i)+\mu_i^2\,(v_1\cdot \bar{h}_i+\phi_i)=(f+\mu_i^2\,v_1)\cdot \bar{h}_i\ \mbox{ in }(0,T)\ \mbox{ a.e. }x\in\Om
\\ \ecart
(v_1\cdot \bar{h}_i+\phi_i)(0,x)=u^0(x)\cdot \bar{h}_i,\ \partial_t(v_1\cdot \bar{h}_i+\phi_i)(0,x)=v^0(x)\cdot\bar{h}_i.
\ea\right.
\eeq
which leads us to
\[
\ba{ll}
\dis (v_1\cdot \bar{h}_i+\phi_i)(t,x) & \dis =\int_0^t{\sin(\mu_i(t-s))\over\mu_i}\,\big(f(s,x)+\mu_i^2\,v_1(s,x)\big)\cdot \bar{h}_i\,ds
\\ \ecart
& \dis +\,\cos(\mu_i t)\,u^0(x)\cdot\bar{h}_i+{\sin(\mu_i t)\over\mu_i}\,v^0(x)\cdot\bar{h}_i.
\ea
\]
Integrating by parts and again using \eqref{CIL1} this yields
\[
\ba{ll}
\dis \phi_i(t,x) & \dis ={\sin(\mu_i t)\over\mu_i}\,\bar{h}_i\cdot v^0(x)+\int_0^t{\sin(\mu_i(t-s))\over\mu_i}\,\bar{h}_i\cdot f(s,x)\,ds
\\ \ecart
& \dis -\int_0^t \cos(\mu_i(t-s))\,\bar{h}_i\cdot\partial_sv_1(s,x)\,ds
\ea
\]
Hence, by summing with respect to $i$ we get that
\[
\ba{ll}
\dis v_2(t,x,y) & \dis =\sum_{i=1}^\infty{\sin(\mu_i t)\over\mu_i}\big(h_i(y)\otimes\bar{h}_i\big)v^0(x)\,ds
+\sum_{i=1}^\infty\int_0^t{\sin(\mu_i(t-s))\over\mu_i}\big(h_i(y)\otimes\bar{h}_i\big)f(s,x)\,ds
\\ \ecart
& \dis -\sum_{i=1}^\infty\int_0^t \cos(\mu_i(t-s))\big(h_i(y)\otimes\bar{h}_i\big)\partial_sv_1(s,x)\,ds
\ea
\]
Finally, defining the kernel
\beq\label{fbK}
\bar{K}(t,y):=\sum_{i=1}^\infty{\sin(\mu_it)\over\mu_i}\,h_i(y)\otimes\bar{h}_i,\quad\mbox{for }(t,y)\in (0,T)\times Y_2,
\eeq
we obtain
\beq\label{u2bK}
\left\{
\ba{l}
\dis v_2(t,x,y)=\bar{K}(t,y)\,v^0(x)+\int_0^t \bar{K}(t-s,y)\,f(s,x)\,ds-\int_0^t\partial_t\bar{K}(t-s,y)\,\partial_s u_1(s,x)\,ds
\\ \ecart
\mbox{a.e. }(t,x,y)\in Q_T\times Y_2.
\ea\right.
\eeq
We have replaced in \eqref{u2bK} the function $v_1$ by the function $u_1$ which are connected by \eqref{v1v2}, since that for a.e. $(t,y)\in (0,T)\times Y_2$ the range of $\bar{K}(t,y)$ is contained in the space spanned by~$b(y)$.
On the other hand, note that using the series expansion \eqref{fbK} and 
\[
\sum_{i=1}^\infty|\bar h_i|^2<\infty,
\]
we can check that
\[
\bar{K}\in L^\infty(0,T;V_2)^3\cap W^{1,\infty}(0,T;H_2)^3\cap W^{2,\infty}(0,T;V_2')^3.
\]
Moreover, since $V_2\subset H^1_0(Y_2)^3$ and the range of $\bar{K}$ is contained in the space spanned by $b$, the kernel satisfies the regularity \eqref{bK}. Formula \eqref{u2bK} also gives an expression of $u_2$, since by \eqref{u1u2u0} and \eqref{v1v2} we have
\beq\label{u2v2b}
u_2=v_2+\left(I-{b\otimes b\over|b|^2}\right)(u^0-u_1).
\eeq
\par
Let us now compute the function $u_3$ in problem \eqref{pb2ss}. We choose $\varphi_1=\varphi_2=0$. We get
$$
\ba{l}\dis \int_{Q_T\times Y_1} {\bf A}_1\big({\bf e}_x(u_1)+{\bf e}_y(u_3)\big):{\bf e}_y(\varphi_3)\,dtdxdy
+\int_{Q_T\times Y_1} (b\times \partial_tu_1)\cdot \varphi_3\,dtdxdy=0.
\ea
$$
Let $w_{jk}$ and $\vartheta_j$, $1\leq j,k\leq 3$, be the vector-valued functions defined by
\beq\label{dewij} \left\{\ba{l}\dis w_{jk}\in H^1_\sharp (Y_1)^3\\ \ecart\dis
\int_{Y_1} {\bf A}_1\big(E_{jk}+{\bf e}_y(w_{jk})\big):{\bf e}_y(\psi)\,dy
=0,\quad\forall\, \psi\in H^1_\sharp (Y_1)^3,
\ea\right.\eeq
where $(E_{jk})_{1\leq j,k\leq 3}$ is the canonical basis in $\RR^{3\times 3}_s$,
\beq\label{devj}
\left\{\ba{l}
\dis \vartheta_j\in H^1_\sharp (Y_1)^3
\\ \ecart
\dis \int_{Y_1}{\bf A}_1{\bf e}_y(\vartheta_j):{\bf e}_y(\psi)\,dy+\int_{Y_1} (b\times e_j)\cdot \psi\,dy=0,
\quad\forall\, \psi\in H^1_\sharp (Y_1)^3.
\ea\right.
\eeq
Then, defining ${\bf W}(y):\RR^{3\times 3}\to \RR^3$ and $V(y)\in \RR^{3\times 3}$ by
\beq\label{defWV}
{\bf W}(y)M:=\sum_{j,k=1}^3 m_{jk}\,w_{jk}(y),\quad V(y)\eta:=\sum_{j=1}^3 \eta_j\,\vartheta_j(y),
\quad \forall\,M\in\RR^{3\times 3},\ \forall\, \eta\in\RR^3,
\eeq
the function $u_3$ is given by
\beq\label{defu3}
u_3(t,x,y)={\bf W}(y)\,{\bf e}_x(u_1)(t,x)+V(y)\,\partial_tu_1(t,x)\ \mbox{ a.e. }(t,x,y)\in Q_T\times Y_1.
\eeq
\subsubsection*{Case where the magnetic has one direction on the boundary of the inclusion}
Assume that $b_{|\partial Y_2}$ has a fixed direction $\xi$ with $|\xi|=1$.
Then, by \eqref{regui} and \eqref{regb} there exists a scalar function $\alpha\in W^{1,\infty}(0,T;L^2(\Om))\times L^\infty(0,T;H^1_0(\Om))$ such that \eqref{defalph} holds.
For any $\beta\in W^{2,\infty}(0,T;L^2(\Om))^3\times W^{1,\infty}(0,T;H^1_0(\Om))^3$ with $\be(0,x)=\be(T,x)=0$, we define 
\beq\label{ph12xi}
\left\{\ba{ll}
\dis \varphi_1(t,x):=\beta(t,x)\,\xi & \mbox{for }(t,x,y)\in Q_T\times Y
\\ \ecart
\dis \varphi_2(t,x,y):=-\left(I-{b(y)\otimes b(y)\over|b(y)|^2}\right)\varphi_1(t,x) & \mbox{for }(t,x,y)\in Q_T\times Y_2
\\ \ecart
\dis \varphi_2(t,x,y):=0 & \mbox{for }(t,x,y)\in Q_T\times Y_1.
\ea\right.
\eeq
Taking $\ph_3=0$ in \eqref{pb2ss} we have
\[
\ba{l}
\dis -\int_{Q_T\times Y}(\partial_tu_1+\partial_tu_2)\cdot (\partial_t\varphi_1+\partial_t\varphi_2)\,dtdxdy
\\ \ecart
\dis +\int_{Q_T\times Y_1} {\bf A}_1\big({\bf e}_x(u_1)+{\bf e}_y(u_3)\big):{\bf e}_x(\varphi_1)\,dtdxdy
+\int_{Q_T\times Y_2}\hskip -8pt {\bf A}_2{\bf e}_y(u_2):{\bf e}_y(\varphi_2)\,dtdxdy\\ \ecart\dis-\int_{Q_T\times Y_1}(b\times u_3)\cdot \partial_t\varphi_1\,dtdxdy=\int_{Q_T\times Y}f\cdot (\varphi_1+\varphi_2)\,dtdxdy.
\ea
\]
Since by \eqref{defalph} $u_1=u^0+\al\,\xi$ and by \eqref{u1u2u0}
\[
u_1+u_2-u^0={b\otimes b\over|b|^2}\,(u_1+u_2-u^0),
\]
by the definitions \eqref{v1v2} of $v_2$ and \eqref{hb} of $\hat{b}$ we also have
\beq\label{uiv2al}
u_1(t,x)+u_2(t,x,y)=u^0(x)+v_2(t,x,y)+\al(t,x)\,\hat{b}(y)\ \mbox{ a.e. }(t,x,y)\in Q_T\times Y_2.
\eeq
Then, using the expressions \eqref{defu3} of $u_3$ and \eqref{ph12xi} of $\ph_1,\ph_2$, and \eqref{bK1} we get
\[
\ba{l}
\dis -\int_{Q_T}\left(|Y_1|+\int_{Y_2}|\hat{b}|^2\,dy\right)\partial_t\al\,\partial_t\be\,dtdx
+\int_{Q_T}\left(\int_{Y_2}v_2\cdot \hat{b}\,dy\right)\partial^2_{tt}\be\,dtdx
\\ \ecart
\dis +\int_{Q_T}{\bf A}_1^*{\bf e}_x(u^0+\al\,\xi):{\bf e}_x(\be\,\xi)\,dtdx
+\int_{Q_T}\partial_t\al\,V_1^*:{\bf e}_x(\be\,\xi)\,dtdx
\\ \ecart
\dis 
+\int_{Q_T\times Y_2}{\bf A}_2{\bf e}_y(\al\,\hat{b}+v_2):{\bf e}_y(\be\,\hat{b})\,dtdxdy
\\ \ecart
\dis-\int_{Q_T}\big({\bf w}^*{\bf e}_x(\al\,\xi)+m^*\,\partial_t \al\big)\,\partial_t\be\,dtdx
=\int_{Q_T}f\cdot\left(|Y_1|\,\cdot\xi+\int_{Y_2}\hat{b}\,dy\right)\be\,dtdx.
\ea
\]
where ${\bf A}_1^*\in \L(\RR^{3\times 3}_s)$, $V_1^*\in\RR^{3\times 3}_s$, ${\bf w}^*:\RR^{3\times 3}\to\RR$, $m^*$ are the homogenized quantities  defined by
\beq\label{A1V1wm*}
\left\{\ba{ll}
\dis {\bf A}_1^*\,E_{jk}:=\int_{Y_1}{\bf A}_1\big(E_{jk}+e_y(w_{jk})\big)\,dy, & 1\leq j,k\leq 3,
\\ \ecart
\dis V_1^*:=\sum_{j=1}^3 \xi_j\int_{Y_1}{\bf e}_y(\vartheta_j)\,dy. &
\\ \ecart
\dis {\bf w}^*E_{jk}:=\xi\cdot\int_{Y_1}b\times w_{jk}\,dy & 1\leq j,k\leq 3
\\ \ecart
\dis m^*:=\xi\cdot\int_{Y_1}b\times(V\xi)\,dy=\sum_{j,k=1}^3\left(\int_{Y_1}{\bf A}_1{\bf e}_y(\vartheta_j):{\bf e}_y(\vartheta_k)\,dy\right)\xi_j\,\xi_k.
\ea\right.
\eeq
This can also be written as
\[
\ba{l}
\dis -\int_{Q_T}\left(|Y_1|+\int_{Y_2}|\hat{b}|^2\,dy\right)\partial_t\al\,\partial_t\be\,dtdx
+\int_{Q_T}\left(\int_{Y_2}v_2\cdot \hat{b}\,dy\right)\partial^2_{tt}\be\,dtdx
\\ \ecart
\dis +\int_{Q_T}{\bf A}_1^*({\bf e}_x(u^0)+\nabla_x\al\odot\xi):(\nabla_x\be\odot\xi)\,dtdx
+\int_{Q_T}\partial_t\al\,V_1^*:(\nabla_x\be\odot\xi)\,dtdx
\\ \ecart
\dis +\int_{Q_T}\left(\int_{Y_2}{\bf A}_2{\bf e}_y(\hat{b}):{\bf e}_y(\hat{b})\,dy\right)\al\,\be\,dtdx
+\int_{Q_T}\left(\int_{Y_2}{\bf A}_2{\bf e}_y(v_2):{\bf e}_y(\hat{b})\,dy\right)\be\,dtdx
\\ \ecart
\dis -\int_{Q_T}\big({\bf w}^*(\nabla_x\al\odot\xi)+m^*\,\partial_t \al\big)\,\partial_t\be\,dtdx
=\int_{Q_T}f\cdot\left(|Y_1|\,\xi+\int_{Y_2}\hat{b}\,dy\right)\be\,dtdx.
\ea
\]
Defining
\beq\label{Malamu*}
\left\{\ba{l}
\dis M^*:=|Y_1|+m^*+\int_{Y_2}|\hat{b}|^2\,dy
\\ \ecart
\dis c^*:=\int_{Y_2}{\bf A}_2{\bf e}_y(\hat{b}):{\bf e}_y(\hat{b})\,dy
\\ \ecart
\dis \la^*\cdot\zeta:={\bf w}^*(\xi\odot\zeta)-V_1^*\xi\cdot\zeta,\ \mbox{ for }\zeta\in\RR^3
\\ \ecart
\dis \mu^*:=|Y_1|\,\xi+\int_{Y_2}\hat{b}\,dy\\ \ecart\dis
A_1^*\zeta:={\bf A}_1^*(\zeta\odot\xi)\,\xi,\ \mbox{ for }\zeta\in\RR^3,
\ea\right.
\eeq
and using the representation \eqref{u2bK} of $v_2$ the previous variational formulation leads us to the following distributional equation
\[
\ba{l}
\dis \partial_{tt}(M^*\al)-\partial_{tt}\left[\int_0^t\left(\int_{Y_2}\partial_t\bar{K}(t-s,y):(\hat{b}(y)\odot\xi)\,dy\right)\partial_s \al(s,x)\,ds\right]
\\ \ecart
\dis +\,\la^*\cdot\nabla_x(\partial_t\al) -\div_x\big(A_1^\ast\nabla_x\alpha\big)+\,c^*\al
-\int_{Y_2}{\bf A}_2{\bf e}_y\left(\int_0^t\partial_s\al(s,x)\,\partial_t\bar{K}(t-s,y)\,\xi\,ds\right):{\bf e}_y(\hat{b})\,dy
\\ \ecart
\dis =-\,\partial_{tt}\left[\int_{Y_2}\bar{K}(t,y):\big(\hat{b}(y)\otimes v^0(x)\big)\,dy\right]
-\int_{Y_2}{\bf A}_2{\bf e}_y\big(\bar{K}(t,y)\,v^0(x)\big):{\bf e}_y(\hat{b})\,dy
\\ \ecart
\dis +\,\mu^*\cdot f-\partial_{tt}\left[\int_0^t\left(\int_{Y_2}\bar{K}(t-s,y)\,f(s,x)\,dy\right)\cdot\hat{b}(y) \,ds\right]
\\ \ecart
\dis -\int_{Y_2}{\bf A}_2{\bf e}_y\left(\int_0^t \bar{K}(t-s,y)f(s,x)\,ds\right):{\bf e}_y(\hat{b})\,dy+\div_x\big({\bf A}_1^*{\bf e}_x(u^0)\xi\big).
\ea
\]
which by the definition \eqref{bK} of the kernel $\bar{K}$ also can be written as
\[
\ba{l}
\dis \partial_{tt}\left[M^*\al-\int_0^t\left(\int_{Y_2}\partial_t\bar{K}(t-s,y):(\xi\odot\xi)\,dy\right)\partial_s \al(s,x)\,ds\right]
\\ \ecart
\dis +\,\la^*\cdot\nabla_x(\partial_t\al)-\div_x\big(A_1^\ast\nabla_x\alpha\big)+\,c^*\al
-\int_{Y_2}{\bf A}_2{\bf e}_y\left(\int_0^t\partial_s\al(s,x)\,\partial_t\bar{K}(t-s,y)\,\xi\,ds\right):{\bf e}_y(\hat{b})\,dy
\\ \ecart
\dis =-\,\partial_{tt}\left[\int_{Y_2}\bar{K}(t,y):\big(\xi\otimes v^0(x)\big)\,dy\right]
-\int_{Y_2}{\bf A}_2{\bf e}_y\big(\bar{K}(t,y)\,v^0(x)\big):{\bf e}_y(\hat{b})\,dy+\mu^*\cdot f
\\ \ecart
\dis -\,\partial_{tt}\left[\int_{Y_2}\left(\int_0^t\bar{K}(t-s,y)\,f(s,x)\,ds\right)\cdot\xi\,dy\right]
-\int_{Y_2}{\bf A}_2{\bf e}_y\left(\int_0^t \bar{K}(t-s,y)\,f(s,x)\,ds\right):{\bf e}_y(\hat{b})\,dy
\\ \ecart
\dis +\,\div_x\big({\bf A}_1^*{\bf e}_x(u^0)\xi\big).
\ea
\]
This provides the homogenized equation \eqref{homeqal} satisfied by $u_1(t,x)=u^0(x)+\al(t,x)\,\xi$.
\subsubsection*{Case where the magnetic has two directions on the boundary of the inclusion}
Finally, assume that $b_{\mid\partial Y_2}$ has two independent directions.
Due to the regularity of $b$ equality \eqref{Condfl1} yields
\[
b(y)\times \partial_tu_1(t,x)=0\ \hbox{ a.e. }(t,x,y)\in Q_T\times \partial Y_2,
\]
which clearly implies \eqref{Casuinm}. Moreover, the proof of formula~\eqref{u2bK1} is quite similar to the proof of~\eqref{u2bK2} in the previous case.
\par\bigskip
It remains to prove the uniqueness of the solution $\al$ to equation \eqref{homeqal}.
To this end, consider a solution $\om\in W^{1,\infty}(0,T;L^2(\Om))\cap L^\infty(0,T;H^1_0(\Om))$ of equation \eqref{homeqal} with nul right-hand side, {\em i.e.}
\[
\left\{\ba{l}
\dis \partial_{tt}\left[M^*\om-\int_0^t\bar{K}_1(t-s)\,\partial_s \om(s,x)\,ds\right]
+\la^*\cdot\nabla_x(\partial_t\om)-\,{\rm div}_x(A_1^*\nabla_x\om)
\\ \ecart
\dis +\,c^*\om-\int_{Y_2}{\bf A}_2{\bf e}_y\left(\int_0^t\partial_s\om(s,x)\,\partial_t\bar{K}(t-s,y)\,\xi\,ds\right):{\bf e}_y(\hat{b})\,dy=0
\ \mbox{ in }Q_T
\\ \ecart
\dis \om(0,\cdot)=0\ \mbox{ in }\Om.
\ea\right.
\]
Then, going back up the former calculations, the functions $z_1$, $z_2$, $z_3$ given respectively from the definitions \eqref{defalph}, \eqref{u2bK2}, \eqref{defu3} of $u_1$, $u_2$, $u_3$, by
\[
\left\{\ba{ll}
\dis z_1(t,x)=\om(t,x)\,\xi, & 
\\ \ecart
\dis z_2(t,x,y)=-\int_0^t\partial_t\bar{K}(t-s,y)\,\partial_s z_1(s,x)\,ds-\left(I-{b(y)\otimes b(y)\over|b(y)|^2}\right)\om(t,x)\,\xi,
\\ \ecart
\dis z_3(t,x,y)={\bf W}(y)\,{\bf e}_x(z_1)(t,x)+V(y)\,\partial_t z_1(t,x),
\\ \ecart
\dis \mbox{a.e. }(t,x,y)\in Q_T\times Y_2,
\ea\right.
\]
are solutions to the variational problem \eqref{pb2ss0} whose solutions are given by \eqref{zi=0}.
Hence, we obtain that $\om(t,x)=0$ a.e. $(t,x)\in Q_T$.
\par
The proof of Theorem~\ref{thm.hom} is now complete.
\subsection{Proof of Proposition~\ref{pro.bK1}}
By \eqref{hj} and the series expansion \eqref{bK} of $\bar{K}$, the scalar function $\bar{k}:=\bar{K}:(\xi\otimes\xi)$ is solution to the equation
\beq\label{eqbk}
\left\{\ba{ll}
\partial^2_{tt}\bar{k}-\div\,(A_2\nabla\bar{k})=0 & \mbox{in }(0,T)\times Y_2
\\ \ecart
\bar{k}(t,\cdot)=0 & \mbox{on }(0,T)\times \partial Y_2
\\ \ecart
\bar{k}(0,\cdot)=0,\ \partial_{t}\bar{k}(0,\cdot)=1 & \mbox{in }Y_2,
\ea\right.
\eeq
where $A_2$ is the definite positive symmetric matrix of $\RR^{3\times 3}$ defined by
\[
A_2\zeta:={\bf A}_2(\zeta\odot\xi)\,\xi,\ \mbox{ for }\zeta\in\RR^3.
\]
By a regularization procedure we may put $1$ as test function in the equation \eqref{eqbk}, which after an integration by parts leads us to the formula
\[
\partial^2_{tt}\left(\int_{Y_2}\bar{k}(t,y)\,dy\right)=\int_{\partial Y_2}A_2\nabla\bar{k}\cdot n\,d\sigma(y).
\]
Then, using the estimate of \cite[Theorem~4.1]{Lio}:
\[
A_2\nabla\bar{k}\cdot n\in L^{\infty}(0,T;L^2(\partial Y_2)),
\]
we get that
\[
\partial^2_{tt}\left(\int_{Y_2}\bar{k}(t,y)\,dy\right)\in L^\infty(0,T).
\]
This combined with definition \eqref{bK1} implies that
\beq\label{bK1bk}
\bar{K}_1(t)=\int_{Y_2}\partial_t\bar{k}(t,y)\,dy\in W^{1,\infty}(0,T).
\eeq
\cqfd


\begin{thebibliography}{20}
\bibitem{BaFiFi} {\sc D.I. Bardzokas, M.L. Filshtinsky \& L.A. Filshtinsky}: {\em Mathematical Methods in Electro-Magneto-Elasticity}, Lecture Notes in Appl. and Compt. Mech. {\bf 32}, Springer-Verlag Berlin 2007, pp. 530.

\bibitem{AGMR}{\sc A.~\'Avila; G.~Griso, B.~Miara, E.~Rohan:} ``Multiscale modeling of elastic waves: theoretical justification and numerical simulation of band gaps", {\em Multiscale Model. Simul.}, {\bf 7} (1) (2008), 1-21.

\bibitem{All}{\sc G.~Allaire:} ``Homogenization and two-scale convergence", {\em SIAM J. Math. Anal.}, {\bf 23} (6) (1992), 1482-1518.

\bibitem{BrCa}{\sc M.~Briane \& J.~Casado-D\'iaz:} ``Increase of mass and nonlocal effects in the homogenization of magneto-elastodynamics problems", arXiv:1806.10998, 2018.

\bibitem{CCMM2}{\sc J.~Casado-D\'{\i}az, J.~Couce-Calvo, F.~Maestre \& J.D.~Mart\'{\i}n-G\'omez}: ``Homogenization and correctors for the wave equation with periodic coefficients", {\em Math. Mod. Meth. Appl. Sci.}, {\bf 24} (2014), 1343-1388.

\bibitem{CoSp}{\sc F.~Colombini \& S.~Spagnolo:} ``On the convergence of solutions of hyperbolic equations", {\em Comm. Partial Differential Equations}, {\bf 3} (1) (1978), 77-103.

\bibitem{FrMu} {\sc G.A. Francfort \& F. Murat}: ``Oscillations and energy densities in the wave equation", {\em Comm. Partial Differential Equations}, {\bf 17} (1992), 1785-1865.
 
\bibitem{Lio}{\sc J.-L. Lions:} {\em Contr\^olabilit\'e exacte, perturbations et stabilisation de syst\`emes distribu\'es}, Tome~1 (French) [Exact controllability, perturbations and stabilization of distributed systems, Vol. 1], Recherches en Math\'ematiques Appliqu\'ees [Research in Applied Mathematics], {\bf 8}, Masson, Paris, 1988, 541 pp.

\bibitem{MiWi}{\sc G.W.~Milton \& J.~Willis:} ``On modifications to NewtonÕs second law and linear continuum elastodynamics", {\em Proc. R. Soc. Lond. Ser. A, Math. Phys. Eng. Sci.}, {\bf 463} (2079) (2007), 855-880.

\bibitem{Ngu}{\sc G.~Nguetseng:} ``A general convergence result for a functional related to the theory of homogenization", {\em SIAM J. Math. Anal.} , {\bf 20} (3) (1989), 608-623.

\bibitem{San}{\sc E.~S\'anchez-Palencia:} {\em Nonhomogeneous media and vibration theory}, Lecture Notes in Physics {\bf 127}, Springer-Verlag, Berlin-New York, 1980, 398 pp.

\bibitem{Sch}{\sc L.~Schwartz:} {\em Mathematics for the physical sciences}, Hermann, Paris, Addison-Wesley Publishing Co., 1966, 358 pp.

\bibitem{Tar1}{\sc L.~Tartar}: ``Homog\'en\'eisation en hydrodynamique", in {\em Singular Perturbation and Boundary Layer Theory}, Lecture Notes in Mathematics, {\bf 597}, Springer, Berlin-Heidelberg 1977, 474-481.

\end{thebibliography}
\end{document}